\documentclass[12pt]{amsart}
\usepackage{epsfig}
\usepackage{amssymb}
\usepackage{graphicx}
\usepackage{color}

\newtheorem{theorem}{Theorem}[section]
\newtheorem{corollary}[theorem]{Corollary}
\newtheorem{lemma}[theorem]{Lemma}
\newtheorem{proposition}[theorem]{Proposition}

\theoremstyle{definition}
\newtheorem{definition}[theorem]{Definition}
\newtheorem{remark}[theorem]{Remark}
\newtheorem{example}[theorem]{Example}

\DeclareMathOperator{\supp}{supp}

\DeclareMathOperator{\rad}{rad}
\DeclareMathOperator{\pos}{pos}
\DeclareMathOperator{\Hom}{Hom}
\DeclareMathOperator{\intt}{int}

\begin{document}

\title{The Toric Hilbert Scheme of a Rank Two Lattice is Smooth and
  Irreducible} 
\author{Diane Maclagan} 
\address{Dept. of Mathematics, Stanford University, Stanford, CA
  94305}
\email{maclagan@math.stanford.edu}
\author{Rekha R. Thomas}
\address{Dept. of Mathematics, University of Washington, Box
  354350, Seattle, WA 98195}
\email{thomas@math.washington.edu}
\thanks{Research partially supported by NSF grant DMS-0100141}
\date{August 4, 2002} 

\begin{abstract} 
  The toric Hilbert scheme of a lattice $\mathcal L \subseteq \mathbb
  Z^n$ is the multigraded Hilbert scheme parameterizing all ideals in
  ${\mathbf k}[x_1,\dots,x_n]$ with Hilbert function value one for
  every $g$ in the grading monoid $G^+=\mathbb N^n/\mathcal L$.  In this
  paper we show that if $\mathcal L$ is two-dimensional, then the
  toric Hilbert scheme of $\mathcal L$ is smooth and irreducible. This
  result is false for lattices of dimension three and higher as the
  toric Hilbert scheme of a rank three lattice can be reducible.
\end{abstract}
\maketitle

\section{Introduction}

The main result of this paper is the following structure theorem.

\begin{theorem} \label{mainthm}
  Let $\mathcal L$ be a two-dimensional lattice contained in $\mathbb
  Z^n$.  Then the toric Hilbert scheme $H_{\mathcal L}$ is smooth and
  irreducible.
\end{theorem}

Consider a sublattice $\mathcal L \subseteq \mathbb Z^n$ of dimension
(rank) $r$ and the abelian group $G = \mathbb Z^n / \mathcal L$.  Let
$S := {\mathbf k}[x_1, \ldots, x_n]$ be the polynomial ring in $n$
variables over an arbitrary infinite field ${\mathbf k}$.  We grade
$S$ by setting $\deg(x_i) = e_i + \mathcal L$ for $i = 1, \ldots, n$
where $e_1, \ldots, e_n$ are the unit vectors of $\mathbb Z^n$. The
set of possible degrees under this grading is $G^+ := \mathbb N^n /
\mathcal L$.

\begin{definition} \label{lgraded}
A homogeneous ideal $I \subseteq S$ is {\bf $\mathcal L$-graded} if
the value of its Hilbert function $\dim_{\mathbf k}((S/I)_g) = 1$ for
all $g \in G^+$.
\end{definition}

The notion of $\mathcal L$-graded ideals extends the notion of
$A$-graded ideals, first introduced in \cite{Arn} and further
developed in \cite{Stu94} and \cite[Chapter 10]{GBCP}. In the
$A$-graded situation, $\mathcal L = \ker_{\mathbb Z}(A):= \{ u \in
\mathbb Z^n \, : \, Au = 0 \}$ where $A$ is an integer matrix, and an
ideal $I \subseteq S$ is $A$-{\em graded} if and only if it is
$\ker_{\mathbb Z}(A)$-graded in the sense of Definition~\ref{lgraded}.
The {\em toric Hilbert scheme} $H_{A}$ \cite{PeSt1}, \cite{Stu94}
parameterizes all $A$-graded ideals for a given $A$. Haiman and
Sturmfels \cite{HaSt} have introduced {\em multi-graded Hilbert
schemes} which provide a uniform setting for many known Hilbert
schemes, including $H_A$. The multi-graded Hilbert scheme is a
quasi-projective scheme which parameterizes all ideals in a polynomial
ring that are homogeneous with respect to grading by a fixed abelian
group and have a fixed Hilbert function.  In the special case where
the Hilbert function takes value one for all elements in the grading
group, Haiman and Sturmfels call the resulting multi-graded Hilbert
scheme a toric Hilbert scheme \cite[Section 5]{HaSt}.

\begin{definition} 
  The {\bf toric Hilbert scheme} of the lattice $\mathcal L$, denoted
  by $H_{\mathcal L}$, is the multi-graded Hilbert scheme that
  parameterizes all $\mathcal L$-graded ideals in $S$.
\end{definition}

In the special case where the two-dimensional lattice $\mathcal L$
equals $\ker_{\mathbb Z}(A)$ ($A$ an integer matrix of corank two),
Theorem~\ref{mainthm} was proved by Gasharov, Peeva, and Stillman
\cite{GaPe}, \cite{PeSt1}.  Another special case already in the
literature is where $\mathcal L$ is a two-dimensional lattice in
$\mathbb Z^2$.  Then $G$ is a finite abelian group, and the lattice
gives an embedding of $G$ into $GL(2)$.  In this case the toric
Hilbert scheme is Nakamura's $G$-Hilbert scheme.  This can be seen by
comparing Reid's functorial description \cite[Section 4.1]{CrawReid}
with that of the toric Hilbert scheme given in \cite{HaSt}.  The fact
that the $G$-Hilbert scheme is smooth and irreducible for abelian
subgroups of $GL(2)$ is due to Kidoh \cite{Kid}.

Thus Theorem \ref{mainthm} can be viewed as providing a common
generalization of both of these results, as well as a more
combinatorial proof of the Gasharov-Peeva-Stillman result.  Although
the $G$-Hilbert scheme is smooth and irreducible for abelian subgroups
of $SL(3)$ \cite{BKR}, \cite{CrawReid} there is no hope for a further
common generalization, as \cite[Theorem 10.4]{GBCP} shows that the
toric Hilbert scheme of a rank three lattice can be reducible.

For an $r$-dimensional lattice $\mathcal L \subseteq \mathbb Z^n$ and
a vector $u \in \mathcal L$, we write $u = u^+ - u^-$ where $u^+, u^-
\in \mathbb N^n$ are defined by setting $(u^+)_i = u_i$ if $u_i>0$,
and $(u^+)_i=0$ otherwise, and $u^- = (-u)^+$.  If $u \in \mathbb
N^n$ we write $x^u$ for the monomial $\prod_{i=1}^n x_i^{u_i}$.  The
{\em lattice ideal} of $\mathcal L$ is the $(n-r)$-dimensional
binomial ideal
$$I_{\mathcal L} = \langle x^{u^+} - x^{u^-} \, : \, u = u^+ - u^- \in
\mathcal L \rangle \subseteq S.$$

If $\mathcal L = \ker_{\mathbb Z} (A)$ for some $A \in \mathbb Z^{d
  \times n}$ of corank $n-d = r$, then $\mathcal L$ is saturated and
  $I_{\mathcal L}$ is the toric ideal of $A$, denoted by $I_A$.  The
  ideal $I_{\mathcal L}$ is a distinguished point on the toric Hilbert
  scheme $H_{\mathcal L}$. It lies on an irreducible component of
  $H_{\mathcal L}$ called the {\em coherent component} \cite[page
  30]{HaSt}.  The algebraic torus $({\mathbf k}^*)^n$ acts on
  $\mathcal L$-graded ideals by scaling variables.  This action
  translates into an action on $H_{\mathcal L}$, and the {\em coherent
  component} of $H_{\mathcal L}$ is the closure of the $({\mathbf
  k}^*)^n$-orbit of $I_{\mathcal L}$.  An ideal $J \in H_{\mathcal L}$
  lies on the coherent component, and is thus called {\em coherent},
  if there is some weight vector $w \in \mathbb Z^n$ and a $\lambda
  \in ({\mathbf k}^*)^n$ such that $J = \lambda \cdot in_w(I_{\mathcal
  L})$, where $in_w(I_{\mathcal L})$ is the {\em initial ideal} of
  $I_{\mathcal L}$ with respect to $w$. Note that $I_{\mathcal L} =
  in_0(I_{\mathcal L})$ and hence is coherent. For an arbitrary
  lattice $\mathcal L$, $H_{\mathcal L}$ could have other components.
  Theorem~\ref{mainthm} asserts that when $\mathcal L$ is
  two-dimensional, the coherent component is the unique component of
  $H_{\mathcal L}$ and that it is smooth.

In Section 2 we establish some general results about $\mathcal
L$-graded ideals for lattices $\mathcal L$ of arbitrary dimension,
while Sections 3 and 4 focus on two-dimensional lattices.  Our proof
of Theorem~\ref{mainthm} is in two parts. In Section~3 we show that
all the monomial ideals on $H_{\mathcal L}$ are coherent
(Theorem~\ref{monocase}). Hence all the fixed points of $H_{\mathcal
L}$ under the action of $({\mathbf k}^*)^n$ lie on the coherent
component. Gr\"obner basis theory implies that every irreducible
component of $H_{\mathcal L}$ contains a monomial ideal.  Thus Theorem
\ref{monocase} implies that every irreducible component of
$H_{\mathcal L}$ intersects the coherent component at a monomial ideal
and so $H_{\mathcal L}$ is connected. In Section~4 we show that the
ideals parameterized by the torus fixed curves between two monomial
ideals in $H_{\mathcal L}$ are also coherent (Theorem~\ref{flipscoh})
which lets us prove that the Zariski tangent space at each monomial
$\mathcal L$-graded ideal is two-dimensional (Lemma~\ref{twoflips}).
As the coherent component is itself two-dimensional, it is therefore
smooth.  Since $H_{\mathcal L}$ is both connected and smooth, it is
irreducible.

%%%%%%%%%%%%%%%%%%%%%%%%%%%%%%%%%%%%%%%%%%%%%%%%%%%%%%%%%%%%%%%%%%%%
\section{General Lattice Lemmas} \label{generallemmas}

Before we restrict to the case where $\dim(\mathcal L)$, the dimension
of $\mathcal L$, equals two, we establish some basic results for
monomial $\mathcal L$-graded ideals when $\dim(\mathcal L) = r \leq
n$.

\begin{definition} \label{graver}
A binomial $x^u - x^v \in I_{\mathcal L}$ is a {\bf Graver binomial}
if there is no other binomial $x^{u'}-x^{v'} \in I_{\mathcal L}$ such
that $x^{u'}$ divides $x^u$ and $x^{v'}$ divides $x^v$. The set of
Graver binomials of $I_{\mathcal L}$ is called the {\bf Graver basis}
of $I_{\mathcal L}$ and is denoted by $Gr_{\mathcal L}$. 
\end{definition}

\begin{example} \label{firsteg}
  For the two-dimensional lattice $\mathcal L \subset \mathbb Z^4$
  generated by $(2,0,-2,-2)$ and $(0,1,1,0)$, we have $Gr_{\mathcal
    L}=\{ x_2x_3-1,x_3^2x_4^2-x_1^2,
  x_1^2x_2^2-x_4^2,x_1^2x_2-x_3x_4^2 \}$.  Note that the first two
  binomials come from the two generators of the lattice.
\end{example}

The Graver basis of $I_{\mathcal L}$ is finite and
\cite[Proposition~5.2]{HaSt} shows that it gives rise to a finite set
of determinantal equations that cut out $H_{\mathcal L}$. For $u \in
\mathbb Z^n$, let $\supp(u) = \{ i \, : \, u_i \neq 0 \}$ be the
support of $u$ and for a monomial $x^u$, we define its support
$\supp(x^u) := \supp(u)$.  If $x^u-x^v \in Gr_{\mathcal L}$, then
$x^u$ and $x^v$ have disjoint supports.  We repeatedly use the
following lemma from \cite{MT} and \cite{PeSt1}.  

\begin{lemma} \cite[Lemma 2.4]{MT} \cite[Lemma 2.2]{PeSt1} \label{mingens}
If $I$ is an $\mathcal L$-graded ideal with $x^u-cx^v$ ($c \in
{\mathbf k}$, possibly zero) in some reduced Gr\"obner basis for $I$,
then $x^u-x^v$ is a Graver binomial. If $c = 0$, then we here assume
that $x^v$ is a monomial not in $I$ of the same $\mathcal L$-degree as
$x^u$.
\end{lemma}

Although in \cite{MT} it is assumed that $\mathcal L=\ker_{\mathbb
Z}(A)$ and $\mathcal L \cap \mathbb N^n = \{0\}$, the proof there is
valid for general lattices.  The same is true for other results quoted
later from \cite{MT}.

Fix a matrix $B \in \mathbb Z^{n \times r}$ whose columns form a
$\mathbb Z$-basis for $\mathcal L$. This implies that $\mathcal L = \{
Bz \, : \, z \in \mathbb Z^r \}$ and the map
\begin{equation} \label{phi}
\phi : \mathbb Z^r \rightarrow \mathcal L \quad \text{given by} \quad
z \mapsto Bz
\end{equation}
is bijective. Let $b_i$ denote the $i$th row of $B$.  The vector
configuration $\mathcal B = \{ b_1, \ldots, b_n \} \subset \mathbb
R^r$ is called the {\em Gale diagram} of $\mathcal L$. For a subset
$\tau \subseteq [n] := \{1, \ldots, n\}$ let $B_{\tau}$ be the
submatrix of $B$ whose rows are indexed by $\tau$ and ${\mathcal
B}_{\tau} := \{b_i \, : \, i \in \tau\}$. Write $\pos({\mathcal
B}_{\tau})$ for the cone $\{ xB_{\tau} \, : \, x \geq 0, x \in \mathbb
R^{|\tau|}\} \subseteq \mathbb R^r$. If $\pos(\mathcal B) = \mathbb
R^r$, $\mathcal B$ is {\em cyclic} and otherwise {\em acyclic}. 
If $r = 2$ and $\mathcal B$ is cyclic, we set $b_{n+1} := b_1$.

\begin{example} \label{bconfig}
  Let $\mathcal L$ be the lattice from Example \ref{firsteg}.  Then
  the Gale diagram of $\mathcal L$ is $\mathcal B = \{ (2,0), (0,1),
  (-2,1), (-2,0)\}$.  This is illustrated in Figure \ref{firstfig}.

\begin{figure}
\resizebox{10cm}{!}{\includegraphics*{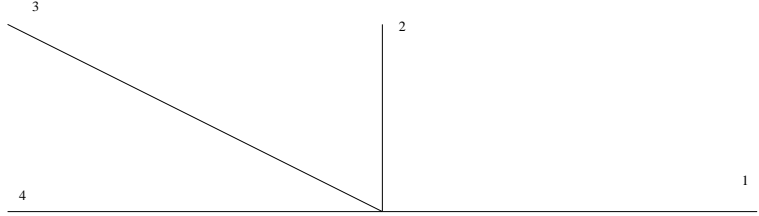}}
\caption{ \label{firstfig} The configuration $\mathcal B$ from Example
\ref{bconfig}}
\end{figure}

\end{example}

\begin{remark}
  We may assume that no $b_i = 0$, and that for $i \neq j$, $b_i \neq
  m b_j$ for any $m \in \mathbb Z_+$. If some $b_i = 0$ then $l_i = 0$
  for all $l \in \mathcal L$ and $x_i$ would not appear in any
  generator of any $\mathcal L$-graded ideal.  This means that there
  would be a bijection between the set of $\mathcal L$-graded ideals
  in $S$ and the set of $\mathcal L'$-graded ideals in ${\mathbf
  k}[x_1,\dots,\hat{x_i}, \dots, x_n]$, where $\mathcal L'$ is the
  projection of $\mathcal L$ which removes the $i$th coordinate.
  Similarly if $b_i=mb_j$ for $m \in \mathbb Z_+$, the map which takes
  every occurrence of $x_j$ to $x_i^mx_j$ would give a bijection
  between $\mathcal L'$-graded ideals and $\mathcal L$-graded ideals.
  These bijections would give rise to an isomorphism between
  $H_{\mathcal L}$ and $H_{\mathcal L'}$, and so it suffices to prove
  Theorem~\ref{mainthm} for the smaller lattice $\mathcal L'$.
\end{remark}

We say that $\pos({\mathcal B}_{\tau})$, or $\tau$, is a $q$-simplex
if $|\tau| = q$ and $\pos(\mathcal B_{\tau})$ is $q$-dimensional. A
{\em triangulation} $T$ of $\mathcal B$ is a collection of
$r$-simplices such that (i) for $\tau, \tau' \in T$, $\pos({\mathcal
B}_{\tau \cap \tau'}) = \pos({\mathcal B}_{\tau}) \cap \pos({\mathcal
B}_{\tau'})$, and (ii) $\pos({\mathcal B}) = \cup_{\tau \in T}
\pos({\mathcal B}_{\tau})$. To be completely accurate, we must
complete $T$ to a simplicial complex on $[n]$ by adding to $T$ all
subsets of $\tau \in T$.  Depending on the context, a simplex $\tau
\in T$ will be either the set $\tau \subseteq [n]$ or the cone
$\pos(\mathcal B_{\tau}) \subseteq \pos(\mathcal B)$.

Recall that all monomial prime  ideals in $S$ are of the form
$P_{\sigma} := \langle x_j \, : \, j \not \in \sigma \rangle$ for a
set $\sigma \subseteq [n]$. The following lemma is a mild extension of
a special case of \cite[Theorem~10.10]{GBCP}.

\begin{lemma} \label{radicallemma}
  Let $I$ be a monomial $\mathcal L$-graded ideal, for a lattice
  $\mathcal L$ with $dim(\mathcal L) = r<n$.  Let $\rad(I)$ be its
  radical and $\rad(I) = \cap_{\sigma \in \Delta(I)} P_{\sigma}$ be
  the unique prime decomposition of $\rad(I)$.  Fix $A = [a_1, \ldots,
  a_n] \in \mathbb Z^{(n-r) \times n}$ such that $AB=0$ and $rank(A) =
  n-r$, and let $\mathcal A = \{a_1, \ldots, a_n \}$.  Then
  $\Delta(I)$ is a triangulation of $\mathcal A$.
\end{lemma}

\begin{proof}
  We first note that since $I_{\mathcal L}$ is $(n-r)$-dimensional,
  the same is true for $I$, so if $\sigma \in \Delta(I)$ we have
  $|\sigma| \leq n-r$.  We will show that $\Delta(I)$ is a
  triangulation of $\mathcal A$, by showing firstly that the
  interiors of any two simplices in $\Delta(I)$ do not intersect in
  this embedding, and secondly that each point in $\pos(\mathcal A)$
  is covered by one of these simplices.  This will also show that
  $|\sigma|=n-r$ for all $\sigma \in \Delta(I)$.
  
  Suppose that $\sigma$ and $\tau$ are two simplices contained in
  simplices in $\Delta(I)$ such that their relative interiors
  intersect.  This means that there is some vector $c \in
  \pos(\mathcal A)$ with $c=\sum_{i \in \sigma} \lambda_i a_i =
  \sum_{j \in \tau} \mu_j a_j$, where we may take $\lambda_i$ and
  $\mu_j$ to be positive integers.  Then $\lambda-\mu \in
  \ker_{\mathbb Z}(A)$, and so there is some multiple $t(\lambda-\mu)
  \in \mathcal L$.  Let $t(\lambda-\mu)=u-v$, where $\supp(u) \cap
  \supp(v) =\emptyset$.  Then since $\supp(u) \subseteq \sigma$, we
  know that $x^u \not \in I$.  Similarly, $\supp(v) \subseteq \tau$,
  so $x^v \not \in I$.  But this means there are two standard
  monomials of $I$ of the same $\mathcal L$-degree, contradicting the
  fact that $I$ is an $\mathcal L$-graded ideal.  This shows that $c$
  does not exist, and so the relative interiors of different simplices
  in $\Delta(I)$ do not intersect.  Note that we made no assumption on
  the dimension of $\sigma$ and $\tau$, so in particular they need not
  be $r$-dimensional.
  
  We now show that each point $c \in \{ u : u = \sum_i n_i a_i, \, n_i
  \in \mathbb Z_+ \}$ is covered by $\pos({\mathcal A}_{\sigma})$ for
  some $\sigma \in \Delta(I)$. Grade the polynomial ring $S$ by
  setting $\deg(x_i)=a_i$.  The $\mathcal L$-grading of $S$ refines
  this grading.  It suffices to show that there is some monomial $x^u$
  of $A$-degree $tc$ for some $t >0$ with $x^u \not \in \rad(I)$.
  This will imply that $\supp(u) \in \Delta(I)$, and so the point $c$
  will be contained in a simplex of $\Delta(I)$.
  
  Consider the ${\mathbf k}$-algebra $U=\oplus_{t\geq 0} S_{tc}$,
  where $S_{tc}$ is the coarsely-graded part of $S$ of degree $tc$.
  We claim that $U$ is a finitely generated algebra, generated by a
  finite number of monomials.  To see this consider the sequence of
  ideals $P_t=\langle x^u : Au=tc \rangle$.  By \cite[Theorem
  1.1]{Mac} only a finite number of the $P_t$ are not contained in
  other ideals of this form.  Let $x^{u_1},\dots, x^{u_s}$ be the
  monomial generators of these inclusion-maximal $P_t$.  If $x^u$ is
  any other monomial of degree $tc$ for some $t>0$, then $x^u$ is
  divisible by one of the $x^{u_i}$.  Since $x^{u-u_i}$ also has
  degree a multiple of $c$, we see that in fact we can write $x^u$ as
  a product of some (possibly repeated) of the $x^{u_i}$, and so
  $x^{u_1},\dots,x^{u_s}$ generate $U$ as an algebra. Let $J=\oplus_{t
  \geq 0} I_{tc}$.  If all of the $x^{u_i}$ lay in the radical of $I$,
  then there would be an $N$ for which $x^{Nu_i} \in I$ for all $i$
  and hence $x^{Nu_i} \in J$ for all $i$.  But this would mean that $U/J$ was a
  finite-dimensional algebra, which contradicts the fact that for all
  $t>0$ there is a standard monomial for $I$ of degree $tc$.  Indeed,
  since the $\mathcal L$-grading refines the $A$-grading, there may
  well be more than one standard monomial of each degree $tc$.  From
  this contradiction we can conclude that there is a generator
  $x^{u_j}$ of $U$ with $x^{u_j} \not \in \rad(I)$, which in turn
  implies that the simplices of $\Delta(I)$ cover $\pos(\mathcal A)$.
\end{proof}

\begin{example} \label{radicaleg}
For the lattice of Example \ref{firsteg}, we have $$B = \left(
\begin{array}{rr}
2 & 0 \\
0 & 1 \\
-2 & 1 \\
-2 & 0 \\
\end{array}\right), \,\, \text{so we can take} \,\, A = \left(
 \begin{array}{llll} 1& -1 & 1 & 0 \\ 0 & -1 & 1 & -1 \\ \end{array}
\right).$$
Note that $\ker_{\mathbb Z}(A)$ is the {\em saturation} of
$\mathcal L$; the lattice generated by the vectors $(1,0,-1,-1)$ and
$(0,1,1,0)$.  With this choice of $A$ we have $\mathcal
A=\{(1,0),(-1,-1), (1,1), (0,-1)\}$.

There are four monomial $\mathcal L$-graded ideals. We list them in
the table below, each with its radical and the prime decomposition of
the radical.\\

$\begin{array}{||c|c|c||}
\hline \text{Monomial} \,\mathcal L \, \text{-graded ideal} &
\text{radical} & \text{prime decomposition of radical}\\ 
\hline
\langle x_2x_3, x_1^2 \rangle & \langle x_2x_3, x_1 \rangle & \langle
x_1,x_2 \rangle \cap \langle x_1,x_3 \rangle\\ 
\langle x_2 x_3, x_1^2x_2,x_3^2x_4^2 \rangle & \langle x_1x_2,
x_2x_3,x_3x_4 \rangle & \langle x_1,x_3 \rangle \cap \langle x_2,x_3
\rangle \cap \langle x_2, x_4 \rangle\\ 
\langle x_2x_3,x_1^2x_2^2,x_3x_4^2 \rangle & \langle x_1x_2,
x_2x_3,x_3x_4 \rangle & \langle x_1,x_3 \rangle \cap \langle x_2,x_3
\rangle \cap \langle x_2, x_4 \rangle\\ 
\langle x_2x_3, x_4^2 \rangle & \langle x_2x_3,x_4 \rangle & \langle
x_2, x_4 \rangle \cap \langle x_3, x_4 \rangle\\ \hline
\end{array}$\\

%%%
%%%There are four monomial $\mathcal L$-graded ideals: $$\langle x_2x_3,
%%%x_1^2 \rangle, \langle x_2 x_3, x_1^2x_2,x_3^2x_4^2 \rangle, \langle
%%%x_2x_3,x_1^2x_2^2,x_3x_4^2 \rangle, \, \text{and} \, \langle x_2x_3,
%%%x_4^2 \rangle.$$
%%%Their radicals are $\langle x_2x_3, x_1 \rangle$ for
%%%the first, $\langle x_1x_2, x_2x_3,x_3x_4 \rangle$ for the next two,
%%%and $\langle x_2x_3,x_4 \rangle$ for the last.  These have primary
%%%decompositions $\langle x_1,x_2 \rangle \cap \langle x_1,x_3 \rangle$,
%%%$\langle x_1,x_3 \rangle \cap \langle x_2,x_3 \rangle \cap \langle
%%%x_2, x_4 \rangle$, and $\langle x_2, x_4 \rangle \cap \langle x_3, x_4
%%%\rangle$ respectively which 
Note that the second and third ideals have the same radical. The three
radicals correspond to the triangulations of $\pos(\mathcal A)$ shown
in Figure \ref{triangspic} in order from left to right.

\begin{figure}
\resizebox{10cm}{!}{\includegraphics*{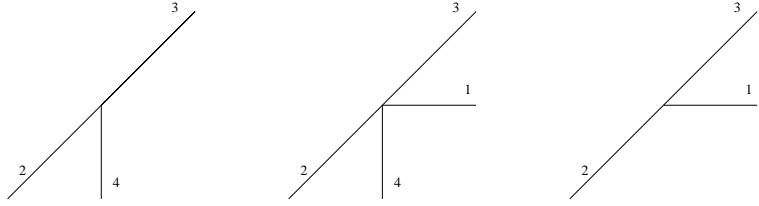}}
\caption{ \label{triangspic} The triangulations of $\mathcal A$
  corresponding to $\mathcal L$-graded ideals.} 
\end{figure} 

\end{example}

\begin{remark} \label{dimension}
  The proof of Lemma \ref{radicallemma} shows that each $\sigma \in
  \Delta(I)$ has cardinality $n-r$. Since the ideals $P_{\sigma}$,
  $\sigma \in \Delta(I)$ are precisely the minimal primes of $I$, this
  shows that $I$ is an equidimensional ideal of dimension $n-r$.  This
  is also true if $r=n$, as then $I$ is a zero-dimensional ideal.
\end{remark}

Let $P_{\sigma}$ be a minimal prime of $I$, so $|\sigma| = n-r$. We
can localize $I$ at $P_{\sigma}$ to get an ideal $I_{\sigma}$. We
identify $I_{\sigma}$ with the projection $\pi_{\sigma}(I)$ where
$\pi_{\sigma}$ is the map:
\begin{eqnarray*}
\pi_{\sigma} : S & \rightarrow & S_{\sigma} := {\mathbf k}[x_i \, : \,
 i \not \in \sigma]\\ x_j & \mapsto & \left\{ \begin{array}{l} x_j
 \quad \text{if} \,\, j \not \in \sigma \\ 1 \quad \text{if} \,\, j
 \in \sigma
\end{array} \right .
\end{eqnarray*}

The ideal $I_{\sigma}$ is also the image under $\pi_{\sigma}$ of the
$P_{\sigma}$-primary component of $I$ from an irredundant primary
decomposition of $I$. Similarly, let $\hat \pi_{\sigma}$ be the map
that projects $u \in \mathbb Z^n$ to the $|\bar \sigma|$-vector
obtained by restricting $u$ to its coordinates indexed by $\bar
\sigma$, where $\bar \sigma = [n] \setminus \sigma$.

\begin{definition}
A homogeneous ideal $J \subseteq S$ is {\bf weakly $\mathcal
    L$-graded} if its Hilbert function $$\dim_{\mathbf k}((S/J)_g) \leq 1
    \,\,\text{for all} \,\, g \in G^+.$$
\end{definition}

We recall (\cite[Lemma 2.6]{MT}) that a monomial ideal $I$ is weakly
  $\mathcal L$-graded if and only if for each $x^u - x^v \in
  Gr_{\mathcal L}$, at least one of $x^u$ or $x^v$ lies in $I$.

\begin{lemma} \label{Isigmafacts}
Let $P_{\sigma}$ be a minimal prime of a monomial $\mathcal L$-graded
ideal $I$.
\begin{enumerate}
\item For each $l \in \mathcal L$ then there is some $i \in
\bar{\sigma}$ for which $l_i\neq 0$.

\item The localized ideal $I_{\sigma} \subset S_{\sigma}$ is weakly
$\mathcal L_{\sigma}$-graded, where $\mathcal L_{\sigma} = \hat
\pi_{\sigma} (\mathcal L) \subseteq \mathbb Z^{|\bar \sigma|}$.  This
implies that $I_{\sigma}$ is an artinian monomial ideal.

\item If $I=in_w(I_{\mathcal L})$, where $w_i =0$ for $i \in \sigma$,
then $I_{\sigma}=in_{\hat \pi_{\sigma}(w)}(I_{\mathcal
L_{\sigma}})$. In particular, $I_{\sigma}$ is coherent.
\end{enumerate}
\end{lemma}

\begin{proof}
\begin{enumerate}
\item Since $P_{\sigma}$ is a minimal prime of $I$, if $x^u \not \in
P_\sigma$ then $x^u \not \in I$. If $l \in \mathcal L$, then $x^{l^+}$
and $x^{l^-}$ have the same $\mathcal L$-degree, so at most one of
them can be a standard monomial of $I$.  This means that at least one
of $x^{l^+}$ and $x^{l^-}$ must lie in $P_{\sigma}$, and so there is
some $i \in \bar{\sigma}$ with $l_i \neq 0$.
  
\item Let $x^u$ and $x^v$ be two monomials in $S_{\sigma}$ with $u-v
\in \mathcal L_{\sigma}$.  Then there exist $u',v'$ with $\supp(u'),
\supp(v') \subseteq \sigma$ and $(u+u')-(v+v') \in \mathcal L$.  But
this means that either $x^{u+u'} \in I$ or $x^{v+v'} \in I$, and thus
one of $x^u$ or $x^v$ must lie in $I_{\sigma}$.  This shows that
$I_{\sigma}$ has at most one standard monomial of each $\mathcal
L_{\sigma}$-degree.

\item By the previous part we know that $I_{\sigma}$ is a weakly
$\mathcal L_{\sigma}$-graded ideal, so it suffices to show that it is
contained in the $\mathcal L_{\sigma}$-graded ideal $in_{\hat
\pi_{\sigma}(w)}(I_{\mathcal L_{\sigma}})$.  Let $x^u$ be a generator
of $I_{\sigma}$.  Then there is a $u' \in \mathbb N^n$ with $\supp(u')
\subseteq \sigma$ and $x^{u+u'} \in I$.  Let $x^v$ be the standard
monomial of $I$ in the same degree as $x^{u+u'}$.  Then $w \cdot
(u+u'-v) >0$, so $\hat \pi_{\sigma}(w) \cdot (u-\hat
\pi_{\sigma}(v))>0$, and so $x^u \in in_{\hat
\pi_{\sigma}(w)}(I_{\mathcal L_{\sigma}})$.

\end{enumerate}
\end{proof}

\begin{remark}
A stronger result than part two of Lemma \ref{Isigmafacts} is true, as
$I_{\sigma}$ is in fact $\mathcal L_{\sigma}$-graded.  The proof is
longer, however, so we prove it in the next section only in the case
where $\mathcal L$ is two-dimensional.
\end{remark}

%%%%%%%%%%%%%%%%%%%%%%%%%%%%%%%%%%%%%%%%%%%%%%%%%%%%%%%%%%%%%%%%%%%
\section{Monomial Ideals}
\label{monomialideals} 

In the rest of this paper we assume that $\dim(\mathcal L) = 2$.  In
this section we show that all monomial ideals in the toric Hilbert
scheme $H_{\mathcal L}$ are coherent.

\begin{theorem} \label{monocase}
Let $I$ be a monomial $\mathcal L$-graded ideal where $\mathcal L$ is a 
two-dimensional sublattice of $\mathbb Z^n$. Then $I=in_w(I_{\mathcal
  L})$ for some $w \in \mathbb Z^n$.
\end{theorem}

The proof of Theorem~\ref{monocase} is established in several steps.
We first show that the localization of a monomial $\mathcal L$-graded
ideal $I$ at a minimal prime $P_{\sigma}$ is coherent in the sense
that it is an initial ideal of $I_{\mathcal L_{\sigma}}$.  One of
these coherent localizations is special in the sense that there is no
other monomial $\mathcal L$-graded ideal with the same
localization. This localization determines $I$. The coherence of
this special localization implies that it is the localization of a
monomial initial ideal of $I_{\mathcal L}$. Thus $I$ is this initial
ideal.

We first recall the result of Lee \cite{Lee} that all triangulations
of $n$ vectors in $\mathbb R^{n-2}$ are {\em regular}.  Recall
\cite[Chapter 8]{GBCP} that a triangulation $\Delta$ of $\{ p_1,
\dots, p_n \} \subseteq \mathbb R^{n-2}$ is regular if there exists a
cost-vector $w \in \mathbb R^n$ such that $\sigma \in \Delta$ if and
only if there exists an $x \in \mathbb Z^{n-2}$ for which $p_i \cdot x
= w_i$ for $i \in \mathbb \sigma$, and $p_i \cdot x < w_i$ for $i \not
\in \sigma$. If such a $w$ exists we denote $\Delta$ as
$\Delta_w$. Applying this definition to triangulations of $\mathcal A$
, we get that $\Delta = \Delta_w$ if and only if $\sigma \in \Delta$
exactly when $wB=(w-xA)B \in \pos(\mathcal B_{\bar \sigma})$.  (Recall
that $AB = 0$.) Thus for a $w \in \mathbb Z^n$ such that $wB \in
\pos(\mathcal B)$, the maximal simplices of the regular triangulation
$\Delta_w$ of $\mathcal A$ are the $(n-2)$-simplices $\sigma \subset
[n]$ such that $wB \in \pos(\mathcal B_{\bar \sigma})$.  This means
that there is a bijection between the regular triangulations of
$\mathcal A$ and the {\em chambers} of $\mathcal B$.

\begin{definition}
  Given a collection $\mathcal P$ of vectors in $\mathbb R^n$, the
  {\bf chamber complex} $\Sigma(\mathcal P)$ of $\mathcal P$ is the
  polyhedral fan obtained by intersecting all the simplices in
  $\mathcal P$.  If $n=2$, then the chamber complex is the collection
  of cones formed by taking the positive hull of adjacent vectors in
  $\mathcal P$.  We identify chambers in the chamber complex with the
  collection of maximal simplices (of $\mathcal P$) which contain
  them.  If $\mathcal L$ is a lattice with generating matrix $B$, we
  denote by $\Sigma(\mathcal L)$ the chamber complex $\Sigma(\mathcal
  B)$.
\end{definition}

\begin{example} For $\mathcal B = \{ (2,0), (0,1), (-2,1), (-2,0)\}$,
the chamber complex is the collection of three cones shown in Figure
\ref{firstfig}.
\end{example}

By Lee's result and Lemma~\ref{radicallemma}, when $\mathcal L$ is
two-dimensional, Lemma \ref{radicallemma} reduces to:

\begin{lemma}
  Let $I$ be a monomial $\mathcal L$-graded ideal where $\dim(\mathcal
  L) = 2$. Then the collection $\{\bar \sigma :
  \sigma \in \Delta(I)\}=\{ \bar \sigma : P_{\sigma} \text{ is a
    minimal prime of } I \}$ is a chamber of $\mathcal B$.
\end{lemma}

\begin{definition} \label{uniquetriang}
  If the monomial $\mathcal L$-graded ideal $I$ maps to the chamber
  $\pos(b_i,b_j)$ then the $(n-2)$-simplex $\sigma = [n] \backslash
  \{i,j\}$ is called the {\bf special simplex} of $I$, and the
  localization $I_{\sigma} \subset S_{\sigma} = {\mathbf k}[x_i,
  x_j]$, the {\bf special localization} of $I$.  If $n=2$ we set
  $I_{\sigma}=I$.
\end{definition}

Note that the special simplex determines the corresponding
triangulation of $\mathcal A$ \cite[Corollary 5.9]{DHSS}. This is
because since $b_i$ and $b_j$ are adjacent in the Gale diagram
$\mathcal B$, $\pos(b_i, b_j)$ does not contain any other chamber of
$\mathcal B$. 

\begin{example} \label{specialsimplexeg}
Let $I$ be the $\mathcal L$-graded ideal $\langle x_2 x_3,
x_1^2x_2,x_3^2x_4^2 \rangle$ from Example~\ref{firsteg}.  We saw in
Example \ref{radicaleg} that $I$ corresponds to the chamber
$\pos(b_2,b_3)$, so the special simplex is $\sigma=\{1,4\}$.  Note
that the corresponding triangulation in Figure \ref{triangspic} is the
only triangulation of $\mathcal A$ containing the simplex
$\pos(a_1,a_4)$. 
\end{example}

From now on, we fix a monomial $\mathcal L$-graded ideal $I$ and let
$\sigma = [n] \setminus \{i,j\}$ always be its special simplex. We now
show that $I_{\sigma}$ is coherent by showing that it is an initial
ideal of $I_{\mathcal L_{\sigma}}$. To prove this it suffices to show
that there exists some $\tilde w \in \mathbb Z^2$ such that $\tilde w
\cdot \hat \pi_{\sigma}(u-v) > 0$ for all $x^u - x^v \in Gr_{\mathcal
L}$ such that $\pi_{\sigma}(x^u) \in I_{\sigma}$ and
$\pi_{\sigma}(x^v) \not \in I_{\sigma}$. This follows from
Lemma~\ref{mingens} because $Gr_{\mathcal L_{\sigma}} \subseteq
\{\pi_{\sigma}(x^u-x^v) : x^u-x^v \in Gr_{\mathcal L}\}$.

\begin{proposition} \label{Isigmacoh}
\begin{enumerate}
\item There exists a cost-vector $w \in \mathbb Z^n$ such that $w
\cdot (u-v)>0$ whenever $x^u-x^v \in Gr_{\mathcal L}$ with
$\pi_{\sigma}(x^v) \not \in I_{\sigma}$.  We can choose $w$ so that
$w_i=0$ for $i \in \sigma$. This implies that $I_{\sigma} = in_{\hat
\pi_{\sigma}(w)}(I_{\mathcal L_{\sigma}})$ and is therefore coherent.
\item The special simplex $\sigma$ of $I$ is also the special simplex
  of the initial ideal $in_w(I_{\mathcal L})$ of $I_{\mathcal L}$, and
  $(in_w(I_{\mathcal L}))_{\sigma}=I_{\sigma}$.
\end{enumerate}
\end{proposition}

\begin{proof}
  (1) Let $x_i^a$ and $x_j^b$ be minimal generators of $I_{\sigma}$
  which exist by Lemma~\ref{Isigmafacts} (2), and let $w \in \mathbb
  N^n$ be the cost-vector with $w_i=b, w_j=a,$ and $w_k=0$ for $k \neq
  i,j$.  Suppose $x^u-x^v \in Gr_{\mathcal L}$, with
  $\pi_{\sigma}(x^v) \not \in I_{\sigma}$.  Since $x^u$ and $x^v$ have
  disjoint support, if $v_i \neq 0$, we must have $u_i=0$.  But then
  we must have $u_j \geq b$, since $\pi_{\sigma}(x^u) \in I_{\sigma}$.
  Since $\pi_{\sigma}(x^v) \not \in I_{\sigma}$, we must have
  $v_i<a$. But now $w \cdot (u-v) =
  -w_iv_i+w_ju_j=-bv_i+au_j>-ba+ab=0$. Similarly, if $v_j \neq 0$ we
  must have $v_j<b$ and $u_i>a$, so $w \cdot (u-v) >0$.  Finally, if
  $v_i=v_j=0$, then $w \cdot (u-v)=bu_i+au_j$. If $u_i=u_j=0$ then
  that would mean that $u-v = l \in \mathcal L$ and $l_i = l_j = 0$
  which contradicts Lemma \ref{Isigmafacts} (1).  So we conclude that
  $w \cdot (u-v) >0$ as required.

\noindent (2) It suffices to show that $\sigma$ is the special simplex
for $in_w(I_{\mathcal L})$, as then Lemma~\ref{Isigmafacts} (3) and
part (1) of this proposition together imply that $(in_w(I_{\mathcal
L}))_{\sigma} = in_{\hat \pi_{\sigma}(w)}(I_{\mathcal L_{\sigma}}) =
I_{\sigma}$. Let $x^v$ be a monomial with $v_i = v_j = 0$. Then for
all $x^u$ such that $x^u - x^v \in I_{\mathcal L}$, either $u_i > 0$
or $u_j > 0$ by Lemma~\ref{Isigmafacts} (1). This implies that
$in_w(x^u-x^v) = x^u$ and hence $x^v \not \in in_w(I_{\mathcal L})$.
This shows that $in_w(I_{\mathcal L}) \subseteq P_{\sigma}$.  Since
$in_w(I_{\mathcal L})$ is $(n-2)$-dimensional, $P_{\sigma}$ must be a
minimal prime of $in_w(I_{\mathcal L})$ and hence $\sigma$ appears in
$\Delta(in_w(I_{\mathcal L}))$. However this implies that $\Delta(I) =
\Delta(in_w(I_{\mathcal L}))$ since $\sigma$ appears in only one
triangulation of $\mathcal A$. Hence $\sigma$ is the special simplex
of $in_w(I_{\mathcal L})$.
\end{proof}

Notice that the proof of Proposition~\ref{Isigmacoh} (1) works for any
localization $I_{\tau}$ of $I$ at a minimal prime $P_{\tau}$ and hence
all these localizations of $I$ are coherent.  Proposition
\ref{Isigmacoh} proves that Theorem \ref{monocase} holds when
$\mathcal L \subseteq \mathbb Z^2$ since in this case, $\mathcal B =
\{b_1, b_2 \}$ and $I =  I_{\emptyset} = I_{\sigma}$.

\begin{example} \label{localeg}
Continuing Example~\ref{specialsimplexeg}, $I_{\sigma}= \langle x_2,
x_3^2 \rangle$, so from the proof of Proposition \ref{Isigmacoh} we
see that $w=(0,2,1,0)$ satisfies $w \cdot (u-v) >0$ whenever $x^u-x^v
\in Gr_{\mathcal L}$ with $\pi(x^v) \not \in I_{\sigma}$.  Every
Graver binomial satisfies this condition on $\pi(x^v)$, and it is easy
to check that they all also satisfy the condition on $w \cdot (u-v)$.
\end{example}

In the rest of this section we fix $w$ to be the vector constructed in
the proof of Proposition~\ref{Isigmacoh} for the special localization
$I_{\sigma}$. The final step in the proof of Theorem~\ref{monocase} is
to show that $I = in_w(I_{\mathcal L})$. This requires understanding
the {\em Gr\"obner fan} of $I_{\mathcal L}$ \cite{BaMo}, \cite{MR}.
This is a fan in $\mathbb R^n$ whose cells (which are open polyhedral
cones) are in bijection with the initial ideals of $I_{\mathcal L}$.
The open cone indexing an initial ideal $J$ of $I_{\mathcal L}$ is the
set of all $p \in \mathbb Z^n$ such that $J = in_{p}(I_{\mathcal L})$.
The closure of this cone is called the {\em Gr\"obner cone} of $J$.
Full dimensional Gr\"obner cones index the monomial initial ideals of
$I_{\mathcal L}$.

The Gr\"obner fan of $I_{\mathcal L}$ can be drawn in $\mathbb R^2$ as
follows. By Lemma~\ref{mingens}, the normal vector to a facet ({\em
  wall}) of a full dimensional Gr\"obner cone is an $l \in \mathcal L$
such that $x^{l^+} - x^{l^-} \in Gr_{\mathcal L}$. Using the injective
map $\phi$ from (\ref{phi}) we can represent $l$ by
$\phi^{-1}(l) \in \mathbb R^2$ and the wall with normal $l$ by the
ray in $\mathbb R^2$ with normal $\phi^{-1}(l)$. We will always mean
this two-dimensional fan when we refer to the Gr\"obner fan of
$I_{\mathcal L}$. If a wall of the fan (ray in $\mathbb R^2$) is
$\pos(g)$ for some $g \in \mathbb R^2$, then we always mean the
primitive clockwise normal vector $g^\perp$ to $g$ by the normal to
this wall.  We now recall that the Gr\"obner fan of $I_{\mathcal L}$
is a refinement of $\Sigma(\mathcal L)$.
Let $\mathbb Z \mathcal B$ be the lattice in $\mathbb Z^2$ generated
by the elements of $\mathcal B$. Recall that $\mathbb Z \mathcal B =
\mathbb Z^2$ if and only if $\mathcal L$ is saturated or equivalently,
if $I_{\mathcal L}$ is a toric ideal.

\begin{definition}
  Let $K$ be a two-dimensional pointed rational polyhedral cone in
  $\mathbb R^2$ and $H$ be its Hilbert basis (i.e., H is a minimal
  generating set for the semigroup $K \cap \mathbb Z \mathcal B$). We
  call the fan obtained by subdividing $K$ by drawing in the rays
  ${\mathbb R}_{\geq 0} \cdot h$ for each $h \in H$ the {\bf Hilbert
  refinement} of $K$.
\end{definition}

By Lemma 3.3.3 in \cite{SST}, the Gr\"obner fan of $I_{\mathcal L}$ is
supported on $\pos(\mathcal B)$. Also, for any $\{r,s\} \subseteq [n]$
and $\tau = [n] \backslash \{r,s\}$ such that $B_{\{r,s\}}$ is
non-singular, the Gr\"obner fan of $I_{\mathcal L_{\tau}}$ is
supported on $\pos(b_r, b_s)$, and Example 3.3.4 in \cite{SST} shows
that the Gr\"obner fan of $I_{\mathcal L_{\tau}}$ is the Hilbert
refinement of $\pos(b_r, b_s)$. Theorem 3.3.8 in \cite{SST} implies
the following.

\begin{lemma} \cite{SST} \label{gbfan}
  The Gr\"obner fan of $I_{\mathcal L}$ is the fan obtained by taking
  the Hilbert refinement of each full dimensional cone (chamber) in
  the chamber complex $\Sigma(\mathcal L)$.
\end{lemma}

A two-dimensional rational cone $K$ is {\em unimodular} if the
primitive integer generators of the two extreme rays of the cone form
a basis for the lattice $\mathbb Z \mathcal B$.

The next two results rely on the geometry of the Hilbert basis of a
two dimensional rational polyhedral cone $K$. It is known that the
Hilbert basis elements of $K$ are precisely the lattice points that
lie on the bounded faces of the polyhedron $K' = \{z \in K \cap
\mathbb Z \mathcal B \, : \, z \neq 0 \}$. See
\cite[Proposition~1.19]{Oda} for instance. For
Corollaries~\ref{coneunimod} and \ref{creepinglemma}, assume that the
Hilbert basis elements of a chamber $\pos(b_r,b_s)$ in
$\Sigma(\mathcal L)$ are $g_0 = b_r, g_1, \ldots, g_t, b_s = g_{t+1}$
in the order they occur on the boundary of $K'$ consisting of its
bounded faces.

\begin{corollary} \label{coneunimod}
  Each full dimensional Gr\"obner cone of $I_{\mathcal L}$ is
  unimodular.
\end{corollary}

\begin{proof}
  If $g_k$ and $g_{k+1}$ are adjacent Hilbert basis elements, then
  there is no element of the lattice $\mathbb Z \mathcal B$ in the
  convex hull of $0$, $g_k$, and $g_{k+1}$ other than the three
  vertices.  Indeed, if a lattice point $g$ existed, then the fact
  that the $g_i$ are vertices of a convex polyhedron means that $g$
  cannot be written as an integral combination of any of the $g_i$,
  which would mean that the Hilbert basis was not complete.  Now
  consider the triangle with vertices $g_k, g_{k+1}$ and $g_k +
  g_{k+1}$. If there was a lattice point $v$ in the interior of this
  triangle, then $v = \lambda g_k + \mu g_{k+1}$ where $0 < \lambda,
  \mu < 1$ and $\lambda + \mu > 1$. This implies that $g_k + g_{k+1} -
  v = g_k(1 - \lambda) + g_{k+1}(1 - \mu)$ lies in the interior of the
  convex hull of $0$, $g_k$, and $g_{k+1}$ since $1-\lambda, 1-\mu
  > 0$ and their sum $1-\lambda + 1 - \mu = 2 - (\lambda + \mu) <
  1$. Since this contradicts the earlier observation, we can conclude
  that there are no lattice points in the interior of the
  parallelogram spanned by $g_k$ and $g_{k+1}$ which is hence a
  fundamental domain of $\mathbb Z \mathcal B$.  Therefore, $\pos(g_k,
  g_{k+1})$ is unimodular.
\end{proof}

\begin{example} The Gr\"obner fan for the lattice of Example
  \ref{firsteg} is shown in Figure \ref{Grobnerpic}.  Note that of the
  three cones in the chamber complex, only the middle one contains
  an extra Hilbert basis element, so we get four Gr\"obner cones.

\begin{figure}
  \resizebox{10cm}{!}{\includegraphics*{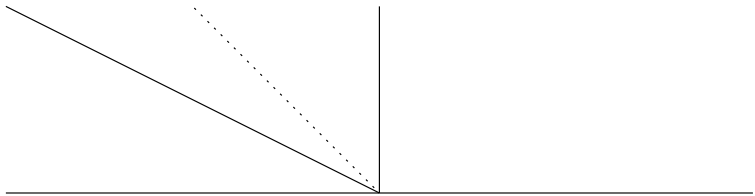}}
\caption{ \label{Grobnerpic} The Gr\"obner fan for the lattice of
  Example \ref{firsteg}.}
\end{figure}
\end{example}

\begin{corollary} \label{creepinglemma}
  Let $\pos(b_r,b_s)$ be a chamber of $\mathcal B$ and $g_0 = b_r,
  g_1, g_2, \ldots, g_t, b_s = g_{t+1}$ be the elements in the Hilbert
  basis of $\pos(b_r,b_s)$ in clockwise order. Then $b_s \cdot
  g_k^\perp < b_s \cdot g_{k+1}^\perp$ for all $k = 1, \ldots, t-1$. 
  Similarly, $b_r \cdot g_k^\perp > b_r \cdot g_{k+1}^\perp$ for all
  $k = 1, \ldots, t-1$. 
\end{corollary}

\begin{proof} Since rotating the cone does not affect the statement, 
  for the first assertion, we may assume without loss of generality
  that $b_r = (0,y)$, $y \in \mathbb Z_+$ and $b_s = (p,q) \in \mathbb
  Z^2$ with $p > 0$.  Let $m_0, m_1, \ldots, m_t$ be the slopes of the
  line segments $[g_0,g_1], [g_1,g_2], \ldots, [g_t,g_{t+1}]$.  Since
  the $g_i$ are vertices of a polyhedron, it now follows that $ m_0
  \leq m_1 \leq \ldots \leq m_t < q/p$. Let $g_k = (a,b)$ and $g_{k+1}
  = (c,d)$.  Then $g_k^\perp = (b,-a)$ and $g_{k+1}^\perp =
  (d,-c)$. Since $m_k = (b-d)/(a-c) < q/p$, we get that $pb-pd <
  qa-qc$ which implies that $b_s \cdot g_k^\perp = pb-qa < pd-qc = b_s
  \cdot g_{k+1}^\perp$. The assertion for $b_r$ is proved similarly.
\end{proof}

We now begin the arguments to prove that $I = in_w(I_{\mathcal L})$.
Recall that $\sigma = [n] \setminus \{i,j\}$ is the special simplex of
$I$. Let $\intt(C)$ denote the interior of a cone $C$.

\begin{lemma}  \label{bklemma}
\begin{enumerate}
\item If $-b_k \in \intt(\pos(b_i,b_j))$, then $k \in \tau$ for all
  $\tau \in \Delta(I)$.
\item If $b_k \cdot \phi^{-1}(l)=0$ for $l \in \mathcal L$, 
  $l_k=0$.
\end{enumerate}
\end{lemma}

\begin{proof}
  (1) Recall that $k \not \in \tau \in \Delta(I)$ if and only if there
  exists some $q$ such that $\pos(b_q,b_k) \supseteq \pos(b_i,b_j)$.
  However, if this was the case, then $-b_k  \in
  \intt(\pos(b_i,b_j))$ would mean  $-b_k \in \intt(\pos(b_k,
  b_q))$ which is not possible. \\
  (2) Let $z=\phi^{-1}(l)$.  Recall that $z \cdot b_i = l_i$ for
  all $i$, and so if $z \cdot b_k = 0$ then $l_k=0$.
\end{proof}

Let $C^{\ast} = \pos(g_1,g_2)$ be the Gr\"obner cone of
$in_w(I_{\mathcal L})$ and $\phi(g_1^\perp)=l_1, \phi(g_2^\perp)=l_2
\in \mathcal L$. Then $C^{\ast}$ is also the Gr\"obner cone of
$(in_w(I_{\mathcal L}))_{\sigma} = in_{\hat \pi( w)}(I_{\mathcal
  L_{\sigma}})$ and hence lies in $\pos(b_i,b_j)$ which is the 
support of the Gr\"obner fan of $I_{\mathcal L_{\sigma}}$. We focus on
the cones $\pos(b_i,b_j)$ and $C^{\ast} = \pos(g_1, g_2)$ in the rest
of this section. Assume they are as in Figure~\ref{cones} (a), where
it could be that $b_i = g_1$ or $b_j = g_2$.

\begin{figure} \label{cones}
\input{cones.pstex_t}
\caption{The cones $\pos(b_i,b_j)$ and $C^{\ast} = \pos(g_1, g_2)$.}
\end{figure}

\begin{lemma} \label{inlemma}
  Suppose $g \in \pos(g_1,g_2)$ such that $g=\lambda g_1 + \mu g_2 \in
  \mathbb Z^2$ with $\lambda, \mu \in \mathbb Z_+$ and
  $\phi(g^\perp)=l \in \mathcal L$. Then $x^{l_1^+}$ divides $x^{l^+}$
  and $x^{l_2^-}$ divides $x^{l^-}$.
\end{lemma}  

\begin{proof}
  We refer to Figure~\ref{cones} (b) for this proof.  If there is some
  index $k$ for which $(l_1)_k$ and $(l_2)_k$ have opposite signs,
  then there must be some vector $l \in \mathcal L \cap \,
  \intt(\pos(l_1,l_2))$ with $l_k=0$.  By Lemma \ref{bklemma} (2),
  $\phi^{-1}(l) \cdot b_k =0$ which implies that $\pm b_k \in
  \intt(C^{\ast})$.  Since $b_k \not \in \intt(C^{\ast})$ for any $k$,
  this can only happen if $-b_k \in \intt(C^{\ast})$.

  Let $U \subseteq [n]$ be the set of indices $k$ for which $-b_k \in
  \intt(C^{\ast})$, where $U$ may be empty.  Then $l,l_1,l_2$ are sign
  compatible in all slots except those indexed by $U$.  If $k \in U$
  then $g_1^\perp \cdot (-b_k) >0$, and so $(l_1)_k<0$.  This means
  that $\supp(l_1^+) \cap U = \emptyset$. Similarly for $k \in U$,
  $g_2^\perp \cdot (-b_k) <0$, so $(l_2)_k >0$, and thus $\supp(l_2^-)
  \cap U = \emptyset$.  Therefore, looking at $l^+-l^- = \lambda
  (l_1^+ - l_1^-) + \mu (l_2^+ - l_2^-) = (\lambda l_1^+ + \mu l_2^+)
  - (\lambda l_1^- + \mu l_2^-) $, we see that $l_1^+ \leq l^+$ and
  $l_2^- \leq l^-$.
\end{proof}

\begin{definition}
  Let $M \subseteq S$ be a weakly ${\mathcal L}$-graded {\em monomial}
  ideal.  The {\bf forced ideal} of $M$ is the ideal generated by all
  monomials $x^u$ for which $x^u-x^v \in Gr_{\mathcal L}$ and $x^v
  \not \in M$.
\end{definition}

\begin{remark} \label{contain}
  If $M$ is a weakly ${\mathcal L}$-graded monomial ideal, then note
  that the forced ideal of $M$ is contained in $M$ as well as in all
  $\mathcal L$-graded monomial ideals contained in $M$. In particular,
  if the forced ideal of $M$ is $\mathcal L$-graded then it is the
  unique $\mathcal L$-graded ideal contained in $M$ since one monomial
  $\mathcal L$-graded ideal cannot be contained in another.
\end{remark}

In Proposition~\ref{Iisforced} we will prove that $I$ is the forced
ideal of $I_{\sigma}$. We now collect together some facts needed for
this proof.

\begin{lemma} \label{Isigmaproperties}
\begin{enumerate}
\item The special localization $I_{\sigma} \subset {\mathbf
k}[x_i,x_j] \subset S$, of the monomial $\mathcal L$-graded ideal $I$,
is weakly $\mathcal L$-graded when considered as an ideal in $S$.
\item The monomials $x^{\hat \pi_{\sigma}(l_1^+)}, x^{\hat
    \pi_{\sigma}(l_2^-)}$ are minimal generators of $I_{\sigma}
  \subset {\mathbf k}[x_i,x_j]$ and the monomials $x^{\hat
    \pi_{\sigma}(l_1^-)}, x^{\hat \pi_{\sigma}(l_2^+)}$ are not in
  $I_{\sigma} \subset {\mathbf k}[x_i,x_j]$. Hence $x^{l_1^+},
  x^{l_2^-}$ lie in the forced ideal of $I_{\sigma} \subset S$.
\item If $g \in \intt(\pos(b_i,b_j))$ and $\phi(g^\perp) = l \in
  Gr_{\mathcal L}$ then $x_i | x^{l^-}$ and $x_j | x^{l^+}$.  If $g
  \not \in \pos(b_i,b_j)$ and $\pos(b_i, b_j)$ lies on the same side
  of $\mathbb R g$ as $g^{\perp}$, then $x^{\hat \pi_{\sigma}(l^-)} =
  1$ which implies that $x^{l^+} \in I_{\sigma} \subset S$.
\end{enumerate}
\end{lemma}

\begin{proof}
\begin{enumerate}
\item If there exists $x^u, x^v \in S$ of the same $\mathcal L$-degree
such that $x^u, x^v \not \in I_{\sigma}$ then $x^u,x^v \not \in I
\subseteq I_{\sigma} \subseteq S$. This cannot happen as $I$ is
$\mathcal L$-graded.

\item Recall that $C^{\ast}$ is the Gr\"obner cone of $in_{\hat
    \pi_{\sigma}(w)}(I_{\mathcal L_{\sigma}})$ with the binomials
  $x^{\hat \pi_{\sigma}(l_1^+)}-x^{\hat \pi_{\sigma}(l_1^-)}$ and
  $x^{\hat \pi_{\sigma}(l_2^-)}-x^{\hat \pi_{\sigma}(l_2^+)}$ defining
  the facets $pos(g_1)$ and $pos(g_2)$ of $C^{\ast}$. Further, the
  initial terms in both binomials are the positive terms $x^{\hat
    \pi_{\sigma}(l_1^+)}$ and $x^{\hat \pi_{\sigma}(l_2^-)}$ which are
  minimal generators of $in_{\hat \pi_{\sigma}(w)}(I_{\mathcal
    L_{\sigma}})$.  It is known that $in_{g_1}(I_{\mathcal
    L_{\sigma}}) = \langle
  x^{\hat \pi_{\sigma}(l_1^+)}-x^{\hat \pi_{\sigma}(l_1^-)} \rangle + 
  \langle x^p : x^p \,\,\text{is a minimal generator of} \,\, in_{\hat
    \pi_{\sigma}(w)}(I_{\mathcal L_{\sigma}}), \, x^p \neq x^{\hat
    \pi_{\sigma}(l_1^+)} \rangle.$ Hence if $x^{\hat
    \pi_{\sigma}(l_1^-)}$ is also in $in_{\hat
    \pi_{\sigma}(w)}(I_{\mathcal L_{\sigma}})$, then we will get that
  $in_{g_1}(I_{\mathcal L_{\sigma}}) \subseteq in_{\hat
    \pi_{\sigma}(w)}(I_{\mathcal L_{\sigma}})$ which is impossible as
  an initial ideal of an ideal cannot be contained in another initial
  ideal of the same ideal.  Therefore, $x^{\hat \pi_{\sigma}(l_1^-)}$
  and similarly, $x^{\hat \pi_{\sigma}(l_2^+)}$ are not in $I_{\sigma}
  = in_{\hat \pi_{\sigma}(w)}(I_{\mathcal L_{\sigma}}) \subset
  {\mathbf k}[x_i,x_j]$. This implies that $x^{l_1^-}, x^{l_2^+}$ are
  not in $I_{\sigma} \subset S$ and hence $x^{l_1^+}, x^{l_2^-}$ lie
  in the forced ideal of $I_{\sigma} \subset S$.

\item Since $l_i = b_i \cdot g^\perp = 0$ if and only if $g$ and $b_i$
are dependent, it follows that $l_i$ changes sign as $g$ is rotated
from one side of $b_i$ to the other. Similarly $l_j$ changes sign as
$g$ is rotated from one side of $b_j$ to the other.  In particular, if
$g \not \in \pos(b_i,b_j)$, then $l_i$ and $l_j$ must have the same
sign, while if $g \in \intt(\pos(b_i,b_j))$, then $l_i$ and $l_j$ must
have opposite signs.
% since otherwise $\pi_{\sigma}(x^{l^+})$ or
%$\pi_{\sigma}(x^{l^-})$ will be the monomial $1$ which cannot be since
%$\pi_{\sigma}(x^{l^+})$ lies in all initial ideals of $I_{\mathcal
%L_{\sigma}}$ corresponding to Gr\"obner cones that lie between $g$ and
%$b_j$ while $\pi_{\sigma}(x^{l^-})$ lies in all initial ideals of
%$I_{\mathcal L_{\sigma}}$ corresponding to Gr\"obner cones that lie
%between $g$ and $b_i$.  
Further, if $g$ lies on the clockwise side of $\pos(b_i,b_j)$ as in
the claim, then $\pi_{\sigma}(x^{l^+})$ lies in all initial ideals
corresponding to Gr\"obner cones in $\pos(b_i,b_j)$ which shows that
$l_i^+, l_j^+ > 0$, as $\pi_{\sigma}(x^{l^-})=1$ which does not lie in
any initial ideal.  Finally, in the case that $g \in \pos(b_i,b_j)$,
the fact that $\pi_{\sigma}(x^{l^+})$ lies in all initial ideals of
$I_{\mathcal L_{\sigma}}$ corresponding to Gr\"obner cones that lie
between $g$ and $b_j$ while $\pi_{\sigma}(x^{l^-})$ lies in all
initial ideals of $I_{\mathcal L_{\sigma}}$ corresponding to Gr\"obner
cones that lie between $g$ and $b_i$ means that $l_i < 0$ and $l_j >
0$.
\end{enumerate}
\end{proof}

\begin{proposition} \label{Iisforced}
  The monomial $\mathcal L$-graded ideal $I$ is the forced ideal of
  $I_{\sigma} \subset S$.
\end{proposition}

\begin{proof}
  Let $J$ be the forced ideal of $I_{\sigma}$. Then $J \subseteq I
  \subseteq I_{\sigma}$. To show that $J = I$, it suffices to show
  that $J$ is weakly $\mathcal L$-graded.  This is true if for every
  $x^u -x^v \in Gr_{\mathcal L}$ either $x^u \in J$ or $x^v \in J$.
  By the definition of $J$, if one of $x^u$ or $x^v$ does not lie in
  $I_{\sigma}$ then the other is in $J$, so we may assume that
  $x^u,x^v \in I_{\sigma}$. Let $g^\perp = \phi^{-1}(u-v)$ and $g \in
  \mathbb Z^2$ be the primitive vector such that $g^\perp$ is the
  clockwise normal to the ray $\pos(g)$.  Since $x^u,x^v \in
  I_{\sigma}$, Lemma~\ref{Isigmaproperties} (3) lets us assume $g \in
  \pos(b_i,b_j)$ (where we use $v-u$ instead of $u-v$ if necessary).
  
  Since $C^{\ast}=\pos(g_1,g_2)$ is unimodular (Corollary
  \ref{coneunimod}), we can write $u-v=\lambda l_1 + \mu l_2$ for
  $\lambda, \mu \in \mathbb Z$. If $\lambda \mu = 0$, then since
  $x^u-x^v \in Gr_{\mathcal L}$ either $u-v = l_1$ or $u-v = l_2$. In
  either case, $x^u$ and $x^v$ do not both belong to $I_{\sigma}$ by
  Lemma~\ref{Isigmaproperties} (2) which contradicts our
  assumption. So $\lambda \mu \neq 0$. The proof now breaks into two
  cases, depending on the sign of $\lambda \mu$.
  
\noindent {\bf Case 1. $\lambda \mu > 0$:} From our assumption that $g
\in \pos(b_i,b_j)$ it follows that $\lambda, \mu >0$.  Then $g \in
\pos(g_1,g_2)$ and hence by Lemma \ref{inlemma}, $x^{l_1^+} | x^u$ and
$x^{l_2^-} | x^v$ which implies that $x^{u}, x^v \in J$ by
Lemma~\ref{Isigmaproperties} (2).
  
\noindent {\bf Case 2. $\lambda \mu <0$:} Then $g \in
\intt(\pos(b_i,b_j)) \setminus C^{\ast}$. We may assume that
$C^{\ast}$ lies on the same side of $\mathbb R g$ as $g^\perp$ (see
Figure~\ref{cones} (c) and (d)) and consider two subcases. 

\noindent {\bf Subcase 2.1:} Suppose $g$ is in the Hilbert basis of
$\pos(b_i,b_j)$. Then by Lemma~\ref{Isigmaproperties} (3), $v_i =
-(b_i \cdot g^\perp) > 0$. Also, $(l_2^-)_i = -(b_i \cdot g_2^\perp)$.
Then Corollary~\ref{creepinglemma} implies that $b_i \cdot g^\perp >
b_i \cdot g_2^\perp$ which implies that $v_i < (l_2^-)_i$. Hence
$x^{\hat \pi_{\sigma}(v)} = x_i^{v_i} \not \in I_{\sigma}$ since
$x_i^{(l_2^-)_i}$ is a minimal generator of $I_{\sigma}$. This means
$x^v \not \in I_{\sigma}$ which is a contradiction.

\noindent {\bf Subcase 2.2:} Suppose $g$ is not in the Hilbert basis
of $\pos(b_i,b_j)$. Since $g$ lies in $\intt(\pos(b_i,b_j)) \backslash
C^{\ast}$, it lies in the interior of some Gr\"obner cone of
$I_{\mathcal L}$ contained in $\pos(b_i,b_j)$ different from
$C^{\ast}$.  Suppose this Gr\"obner cone is $\pos(g_3,g_4)$ with
$\phi(g_3^\perp) = l_3 \in Gr_{\mathcal L}$ and $\phi(g_4^\perp) = l_4
\in Gr_{\mathcal L}$ as in Figure~\ref{cones} (d). By
Lemma~\ref{inlemma}, $x^{l_3^+} | x^u$. Since $g_3$ lies in the
Hilbert basis of $\pos(b_i,b_j)$, the arguments in the previous
subcase show that $x^{l_3^-} \not \in I_{\sigma}$ which implies that
$x^{l_3^+} \in J$ and hence $x^u \in J$.
\end{proof}

\begin{example}
In Example \ref{specialsimplexeg} we saw that for $I=\langle x_2 x_3,
x_1^2x_2,x_3^2x_4^2 \rangle$ we have $\sigma=\{1,4\}$, and
$I_{\sigma}=\langle x_2, x_3^2 \rangle$.  Then in the notation of the
preceding proofs, $b_i=b_3$, $b_j=b_2$, $g_1=(-1,1)$, and $g_2=b_2$.
We now check that the forced ideal of $I_{\sigma}$ is weakly $\mathcal
L$-graded, by checking that each binomial in the Graver basis has one
of its monomials lying in $I_{\sigma}$.  The Graver binomial
$x_2x_3-1$ is covered by part 3 of Lemma \ref{Isigmaproperties}, as $1
\not \in I_{\sigma}$.  The binomial $x_3^2x_4^2-x_1^2$ comes from
$\phi(b_2^{\perp})$, while $x_1^2x_2-x_3x_4^2$ comes from
$\phi(g_1^{\perp})$, so they are both dealt with by part 2 of Lemma
\ref{Isigmaproperties}.  Note that $x_1^2$ and $x_3x_4^2$ do not lie
in $I_{\sigma}$.  Finally, $x_1^2x_2^2-x_4^2$ comes from
$\phi(b_3^{\perp})$, so is covered by subcase 2.1. of
Proposition~\ref{Iisforced}. Note that $x_4^2 \not \in I_{\sigma}$.
\end{example}

The intersection of the minimal primary components of a monomial ideal
$M$ is called $Top(M)$ in \cite{SST}. Theorem~3.3.6 in \cite{SST} and
Theorem 4.4 in \cite{HT1} then follow from the following corollary of
Proposition~\ref{Iisforced}.

\begin{corollary} If $dim(\mathcal L) = 2$ and $I$ and $I'$ are two
  distinct monomial $\mathcal L$-graded ideals, then $Top(I) \neq
  Top(I')$.
\end{corollary}

\begin{proof}
If $I \neq I'$, then Proposition \ref{Iisforced} implies that
$I_{\sigma} \neq I'_{\sigma'}$, where $\sigma$ and $\sigma'$ are the
special simplices of the two ideals.  Since the minimal primary
components are uniquely determined from $Top(I)$, this implies that
$Top(I) \neq Top(I')$.
\end{proof}

Propositions \ref{Isigmacoh} and \ref{Iisforced} combine to prove
Theorem \ref{monocase}.

\begin{proof}[Proof of Theorem \ref{monocase}]
  Proposition \ref{Iisforced} says that $I$ is generated by monomials
  $x^u$ for which there is some $x^v \not \in I_{\sigma}$ with
  $x^u-x^v \in Gr_{\mathcal L}$.  Proposition \ref{Isigmacoh} says
  that there is a cost-vector $w$ for which $w \cdot (u-v) > 0$ for
  all such Graver binomials $x^u-x^v$.  This implies that $I \subseteq
  in_w(I_{\mathcal L})$.  Since we cannot have a proper inclusion of
  monomial $\mathcal L$-graded ideals, we conclude that
  $I=in_w(I_{\mathcal L})$.
\end{proof}

%%%%%%%%%%%%%%%%%%%%%%%%%%%%%%%%%%%%%%%%%%%%%%%%%%%%%%%%%%%%%%%%%%%%%
\section{The General Case}

In this section we finish proving Theorem \ref{mainthm}.  We first
recall from \cite{MT} the notion of a {\em flip}.

\begin{definition}
  Let $M_1$ be a monomial $\mathcal L$-graded ideal with minimal
  generator $x^u$.  Let $x^v$ be the standard monomial for $M_1$ of
  the same degree as $x^u$.  The {\em wall ideal} $W$ associated to
  $M_1$ and $x^u$ is the ideal obtained by replacing the generator
  $x^u$ of $M_1$ by the binomial $x^u-x^v$.  We say the binomial
  $x^u-x^v$ is {\em flippable} if $in_{\prec}(W)=M_1$ for any term
  order with $x^v \prec x^u$.  If there is some other term order
  $\prec'$ with $x^u \prec' x^v$ then we set $M_2$ to be the $\mathcal
  L$-graded monomial ideal $in_{\prec'}(W)$.  We say that $M_2$ and
  $M_1$ differ by a flip over the {\em true flip} $x^u-x^v$.  If there
  is no such term order $\prec'$ we say that $x^u-x^v$ is a {\em fake
  flip}.
\end{definition}

\begin{remark} \label{flipstuff}
\begin{enumerate}
\item In \cite{MT} the assumption that $\ker(A) \cap \mathbb N^n =
  \{0\}$ (the {\em positively-graded} assumption) avoided the
  possibility of fake flips.
  
\item For a true flip, we recall Definition 2.7 and Lemma 2.9 of
  \cite{MT} which says that $M_2=\langle x^{\alpha} \in M_1 : \alpha \neq u \,
  \text{and} \, \exists \, x^{\beta} \not \in M_1 \text{ with }
  x^{\alpha}-x^{\beta} \in G r_{\mathcal L} \rangle + \langle x^v
  \rangle.$

\end{enumerate}
\end{remark}

\begin{example} Let $\mathcal L$ be the lattice of Example
\ref{firsteg}.  For the $\mathcal L$-graded ideal $\langle x_2x_3,
x_1^2 \rangle$ the binomial $x_2x_3-1$ is a fake flip, and the
binomial $x_1^2-x_3^2x_4^2$ is a true flip.  Flipping over the true
flip, we get the other initial ideal of $\langle x_1^2-x_3^2x_4^2,
x_2x_3 \rangle$, which is $ \langle x_2 x_3, x_1^2x_2,x_3^2x_4^2 \rangle$.
This has two true flips: $x_1^2-x_3^2x_4^2$ and $x_1^2x_2-x_3x_4^2$.
Flipping over the second true flip gives the ideal $\langle
x_2x_3,x_1^2x_2^2,x_3x_4^2 \rangle$, which again has two true flips:
$x_1^2x_2-x_3x_4^2$ and $x_1^2x_2^2-x_4^2$.  Flipping over
$x_1^2x_2^2-x_4^2$ we get the last $\mathcal L$-graded ideal $\langle
x_2x_3, x_4^2 \rangle$.  For this ideal $x_1^2x_2^2-x_4^2$ is a true
flip, and $x_2x_3-1$ is a fake flip.
\end{example}

 We first prove that every wall ideal is coherent.
If $x^a$ is a monomial, we write $x^{\supp(a)}$ for the monomial
$\prod_{i \in \supp(a)} x_i$.

\begin{lemma} \label{flipscoh}
Let $\mathcal L$ be a two-dimensional lattice contained in $\mathbb
Z^n$, and let $I$ and $J$ be two monomial initials ideals of
$I_{\mathcal L}$ which differ by a flip.  Then the wall ideal $W$ of
$I$ and $J$ is coherent.
\end{lemma}

\begin{proof}
 Let the flip be over the binomial $x^u-x^v$, with $x^u \in I
\setminus J$, and $x^v \in J \setminus I$.  Choose cost-vectors $w_0$
and $w_1$ for which $in_{w_0}(I_{\mathcal L})=I$, and
$in_{w_1}(I_{\mathcal L})=J$.  Let $w'$ be the composite cost-vector
$(w_0 \cdot (u-v))w_1-(w_1 \cdot (u-v))w_0$.  This is chosen so that
$w' \cdot (u-v)=0$.  Pick a monomial generator $x^{\alpha}$ of $W$.
If the standard monomial of the same degree as $x^{\alpha}$ is the
same monomial $x^c$ for both $I$ and $J$, then $w_i \cdot (\alpha-c)
>0$ for $i=0,1$.  This means that $w' \cdot (\alpha-c)>0$, and so
$x^{\alpha} \in in_{w'}(I_{\mathcal L})$.  If this is the case for all
such minimal monomial generators $x^{\alpha}$, then $W \subseteq
in_{w'}(I_{\mathcal L})$.  As we cannot have a proper inclusion of
$\mathcal L$-graded ideals, we conclude that in this case we have
$W=in_{w'}(I_{\mathcal L})$.

Suppose on the contrary that there is some minimal generator
$x^{\alpha}$ of $W$ for which the standard monomials for $I$ and $J$
are different monomials $x^c$ and $x^d$ respectively.  The rest of the
proof deals with this case.  Let $x^e-x^f$ be a Graver binomial with
$x^e$ dividing $x^c$ and $x^f$ dividing $x^d$.  Then we must have $x^e
\in J \setminus I$ and $x^f \in I \setminus J$.  The second
observation in Remark \ref{flipstuff} now implies that $e=v$ and
$f=u$.  So we can write $c=kv+\gamma$, and $d=ku+\gamma$, where
$\gamma=\gamma_u+\gamma_v+\gamma'$ with $\supp(\gamma_u) \subseteq
\supp(u)$, $\supp(\gamma_v) \subseteq \supp(v)$, and $\supp(\gamma')
\cap (\supp(u) \cup \supp(v)) = \emptyset$.  Also note that by Lemma
\ref{mingens} $x^{\alpha}-x^c$ and $x^{\alpha}-x^d$ are both Graver
basis elements, so $\supp(\alpha) \cap \supp(c) = \supp(\alpha) \cap
\supp(d) = \emptyset$.

We next claim that $x^{u+v+\gamma} \not \in \rad(W)$.  If
$x^{u+v+\gamma} \in \rad(W)$, there is some $l \in \mathbb N$ for
which $x^{l(u+v+\gamma)} \in W$.  We may assume that $l \geq k$.
Choose a monomial $x^g$ with $\supp(g) \subseteq \supp(u)$ so that
$x^{\gamma_u+g}=x^{pu}$ for some integer $p$.  Since $x^u-x^v \in W$,
and $x^{l(u+v+\gamma)+lg} \in W$, it follows that
$x^{l((2+p)v+\gamma_v+\gamma')} \in W$, and thus
$x^{l((2+p)v+\gamma_v+\gamma')} \in I$.  Note that
$\supp(l((2+p)v+\gamma_v+\gamma'))=\supp(v) \cup \supp(\gamma')$.

Let $x^s$ be a monomial not lying in $\rad(I)$.  Since $x^{\alpha} \in
W$, and thus $x^{\alpha} \in I$, we cannot have $\supp(\alpha)
\subseteq \supp(s)$.  Similarly we cannot have $\supp(u) \subseteq
\supp(s)$ or $\supp(v) \cup \supp(\gamma') \subseteq \supp(s)$.  Since
these three sets are pairwise disjoint, we conclude that $|\supp(s)|
\leq n-3$.  But since $s$ was arbitrary this contradicts the
observation of Remark \ref{dimension} that $I$ is a
$(n-2)$-dimensional ideal.  From this we conclude that $x^{u+v+\gamma}
\not \in \rad(W)$.

Since $\alpha-c$ and $\alpha-d$ are not linearly dependent, they must
be a basis for  the real span of $\mathcal L$, so every vector $l \in \mathcal L$
can be written as a (possibly non-integral) linear combination of
$\alpha-c$ and $\alpha-d$.

Let $w$ be the cost-vector with $w_i=1$ if $i \in \supp(\alpha)$, and
$w_i=0$ otherwise.  We will show that $W=in_w(I_{\mathcal L})$.  Let
$F$ be the generating set for $in_w(I_{\mathcal L})$ obtained from the
corresponding reduced Gr\"obner basis for $I_{\mathcal L}$.  If $f \in
F$ then $f$ is either a binomial or a monomial.  If $f$ is a binomial,
then $f$ comes from a lattice vector
$l=\lambda(\alpha-c)+\mu(\alpha-d)$ with $w \cdot l=0$.  This means
that $\lambda=-\mu$, and so $l=ku-kv$.  Since $f$ must be a Graver
binomial, we know that $k=1$, and thus $f=x^u-x^v \in W$.  Now
consider the case that $f$ is a monomial $x^{\beta}$, and let
$x^{\delta}$ be a standard monomial of $in_w(I_{\mathcal L})$ with
$x^{\beta}-x^{\delta} \in Gr_{\mathcal L}$.  Since $\beta-\delta \in
\mathcal L$, it equals $\lambda(\alpha-c) + \mu(\alpha-d)$ for some
constants $\lambda$ and $\mu$, so
$\beta-\delta=(\lambda+\mu)\alpha-k\lambda v -k\mu u
-(\lambda+\mu)\gamma$.  Because $x^{\beta}$ is the leading term of
$x^{\beta}-x^{\delta}$ with respect to the cost-vector $w$, we know
that $\lambda+\mu>0$.  It follows that $\supp(\delta) \subseteq
\supp(u) \cup \supp(v) \cup \supp(\gamma)$, so $x^{\delta} \not \in
W$.  But this means that there is some constant $\rho$ for which
$x^{\beta}-\rho x^{\delta} \in W$, and so $x^{N\beta}-\rho^Nx^{N\delta} \in
W$ for all $N \in \mathbb N$.  We can choose $N$ sufficiently large so
that $N(\lambda+\mu)>1$, so $x^{\alpha}$ divides $x^{\beta}$.  But
$x^{\alpha}$ was chosen to be in $W$, so if $\rho \neq 0$ this implies
that $x^{N\delta} \in W$, and so $x^{\delta} \in \rad(W)$.  From this
contradiction we conclude that $\rho=0$, and so $x^{\beta} \in W$.  This
shows that every minimal generator of the $\mathcal L$-graded ideal
$in_w(I_{\mathcal L})$ lies in $W$, and so since $W$ is itself
$\mathcal L$-graded, we must have $W=in_w(I_{\mathcal L})$.
\end{proof}

\begin{corollary} \label{twoflips}
Let $\mathcal L$ be a two-dimensional lattice.  Then every monomial
$\mathcal L$-graded ideal has exactly two flips.
\end{corollary}

\begin{proof}
Let $I$ be a monomial $\mathcal L$-graded ideal.  By Theorem
\ref{monocase} we know that $I$ is an initial ideal of $I_{\mathcal
L}$, and thus corresponds to a cone in the Gr\"obner fan of
$I_{\mathcal L}$.  Now Lemma \ref{flipscoh} says that any true flip of
$I$ connects $I$ to an initial ideal $J$ that corresponds to an
adjacent Gr\"obner cone.  There are thus exactly two true flips of
$I$, unless $I$ corresponds to a cone which is adjacent to the
boundary of the Gr\"obner fan of $I_{\mathcal L}$, in which case $I$
has only one true flip.

It thus suffices to show that $I$ has no fake flips unless it
corresponds to a cone adjacent to the boundary of the Gr\"obner fan,
in which case $I$ has exactly one fake flip.

We first recall that if $x^u-x^v$ is a flip of $I$ then $x^u-x^v \in
Gr_{\mathcal L}$.  This means that $x^u$ and $x^v$ have disjoint
supports, and so if $\supp(v) \neq \emptyset$ there is some term order
$\prec$ with $x^u \prec x^v$.  Thus the only way for $x^u-x^v$ to be a
fake flip is to have $x^v=1$.

Suppose $x^u-1$ is a fake flip for $I$.  If $x^{u'}$ were another
minimal generator of $I$ with $\supp(u) \cap \supp(u') \neq
\emptyset$, then the $S$-pair $S(x^u-1,x^{u'})$ would divide $x^{u'}$,
contradicting the flippability of $x^u-1$.  On the other hand, if
$x^{u'}$ is a minimal generator of $I$ with $\supp(u) \cap \supp(u') =
\emptyset$, then $S(x^u-1,x^{u'})=x^{u'} \in I$.  So we see that a
necessary and sufficient condition for a minimal generator $x^u$ of
$I$ of the same $\mathcal L$-degree as $1$ to give rise to a fake
flip $x^u-1$ is for $\supp(u) \cap \supp(u')=\emptyset$ for all other
minimal generators $x^{u'}$ of $I$.

We note that for any vector  $x^u-1 \in I_{\mathcal L}$ (such as a fake flip)  there is a vector $v \in
\mathbb Z^2$ for which $u=\phi(v)$, and $b_i \cdot v \geq 0$ for all
$i$.  Note that we can have $b_i \cdot v=0$ for at most two values of
$i$, and if $b_i \cdot v = b_j \cdot v =0$, then $b_i$ is a negative
multiple of $b_j$.

Let the Gr\"obner cone of $I$, $\pos(g_1,g_2)$, be contained in the
secondary cone $\pos(b_i,b_j)$.  Let $g_k^{\perp}$ be the clockwise normal
to $g_k$ for $k=1,2$, and let $l_k=\phi(g_k^{\perp}) \in \mathcal L$.  We
assume that $g_1$ and $g_2$ are in clockwise order.

We first consider the case where $I$ is a monomial $\mathcal L$-graded
ideal with two true flips.  Then $x^{l_1^+}, x^{l_2^-} \in I$.  Note
that $(l_1)_j=g_1^{\perp} \cdot b_j >0$, and $(l_2)_i = g_2^{\perp}
\cdot b_i <0$.  Hence $i,j \not \in \supp(u)$. This means that if
$x^u-1$ is a fake flip with $u=\phi(v)$, then $v \cdot b_j=v \cdot
b_i=0$.  Since $\pos(b_i,b_j)$ is a pointed cone this is not possible
for $v \neq 0$, so we conclude that $I$ has no fake flips.

We now consider the case where $I$ is a monomial $\mathcal L$-graded
ideal with at most one true flip, so $I$ corresponds to a cone on the
boundary of the Gr\"obner fan of $I_{\mathcal L}$.  Without loss of
generality we may assume that this boundary is the counter-clockwise outer
wall $b_i$.  We first note that the fact that the Gr\"obner fan has a
boundary implies that it is pointed.  If the fan is not pointed, then
for every vector $x \in \mathbb Z^2$ there is some $b_k$ with $x \cdot
b_k <0$, and so $\mathcal L \cap \mathbb N^n = \{ 0 \}$.  This means
that the ideal $I_{\mathcal L}$ is positively graded, which in turn
means that the Gr\"obner region of the ideal is all of $\mathbb R^n$,
so the fan does not have a boundary.  The fact that the Gr\"obner fan
is pointed means that $b_i^{\perp} \cdot b_j \geq 0$ for all $i$, and
so if $u=\phi(b_i^{\perp})$, then $u \geq 0$, and we will show that
$x^u-1$ is a fake flip.

We know that the Gr\"obner cone for $I$ is $\pos(b_i,g_2)$ for some
vector $g_2$, and $x^{l_2^-} \in I$.  Note that $l_2^-=ae_i$ for some
$a \in \mathbb N$.  We now show that this $x_i^a$ is the only
generator of $I$ divisible by $x_i$.  Let $l' \in \mathcal L$, and
write $l'=\phi(g^{\perp})$ for a vector $g^{\perp} \in \mathbb Z^2$.
Write $g$ for the primitive vector in $\mathbb Z^2$ whose clockwise
normal is $g^{\perp}$.  By Corollary \ref{coneunimod} we know that
$\pos(b_i,g_2)$ is a unimodular cone, so $g=\lambda b_i + \mu g_2$
with $\lambda, \mu \in \mathbb Z$.  Since $\phi(b_i^{\perp})_i=0$, we
know that
%%$l'_i=\mu \phi(g_2^{\perp})_i = \mu a$.  Since $|\mu a| >a$
%%for $|\mu| \neq 1$, we see that one of $x^{l'^-}$ and $x^{l'^-}$ is
%%in
$l'_i=\mu \phi(g_2^{\perp})_i = -\mu a$.  Since $|\mu a| >a$ for
$|\mu| \neq 1$, we see that $x^{{l'}^{-}}$ is in $I$ but not a minimal
generator, so $l'$ does not give rise to a minimal generator of $I$
unless $l'=l_1$ or $l_2$.  So the two generators of $I$ are
$x^{l_1^+}=x^u$ and $x^{l_2^-}=x_i^a$.  Since $u_i=b_i^{\perp} \cdot
b_i =0$, it follows from the above characterization of fake flips that
$x^u-1$ is a fake flip.  The only other possible flip come from
$x_i^a$.  If $I$ has any true flips, this must come from $x_i^a$,
giving $I$ a total of two flips.  Otherwise $I$ is the only monomial
$\mathcal L$-graded ideal, and $g_2$ must be the other boundary of the
Gr\"obner fan, so $x_i^a$ gives rise to a second fake flip.
\end{proof}

We recall the following theorem which is Corollary 5.2 of \cite{PeSt1}.

\begin{theorem} \label{tangentflips}
Let $I$ be a monomial $\mathcal L$-graded ideal.  Then the number of
flips of $I$ is equal to $\dim_{\mathbf k} (\Hom_S(I,S/I))_0$.
\end{theorem}

We note that in \cite{PeSt1} the assumption was made that $\mathcal L$
was saturated and positively graded.  The proof goes through word for
word in the case of general lattices.

We also recall the following theorem which is Proposition 1.6 of
\cite{HaSt} applied to $H_{\mathcal L}$.

\begin{theorem} \label{tangentspacedim}
The Zariski tangent space to the toric Hilbert scheme $H_{\mathcal L}$
at an ideal $I$ is canonically isomorphic to $(Hom_S(I,S/I))_0$.
\end{theorem}

We now complete the proof of Theorem \ref{mainthm}. 

\begin{proof}[Proof of Theorem \ref{mainthm}]
There is an action of the torus $({\mathbf k}^*)^n$ on $H_{\mathcal L}$
 given by scaling the variables occurring in an ideal $I$.  Since the
 singular locus of $H_{\mathcal L}$ must be fixed under this torus
 action, and the monomial $\mathcal L$-graded ideals are the
 torus-fixed points, to show that the scheme is smooth we need only
 show that it is smooth at each monomial ideal.  Now Corollary
 \ref{twoflips} and Theorem \ref{tangentflips} imply that for a
 monomial $\mathcal L$-graded ideal we have $\dim_{\mathbf k} \Hom(I,S/I)=2$,
 which Theorem \ref{tangentspacedim} says is the dimension of the
 tangent space of $H_{\mathcal L}$ at $I$.    

Pick any two linearly independent vectors $b_i,b_j \in \mathcal B$.
 Let $\lambda(a,b) \in ({\mathbf k}^*)^n$ have $(\lambda(a,b))_l=1$ for $l \neq
 i,j$, $(\lambda(a,b))_i=a$, and $(\lambda(a,b))_j=b$.  There exist
 $l,l' \in \mathcal L$ with $l_i=0, l_j \neq 0$ and $l'_i \neq 0,
 l'_j=0$, so considering the action of $\lambda(a,b)$ on the binomials
 $x^{l^+}-x^{l^-}$ and $x^{{l'}^+}-x^{{l'}^-}$ we see that the map of
 $\mathbb ({\mathbf k}^*)^2$ to the underlying reduced scheme of $H_{\mathcal
 L}$ given by $(a,b) \mapsto \lambda(a,b) I_{\mathcal L}$ is
 injective, so the coherent component of $H_{\mathcal L}$ is at least
 two-dimensional.  Since we showed above that the dimension of the
 tangent space to $H_{\mathcal L}$ at each monomial ideal is two, this
 shows that $H_{\mathcal L}$ is smooth at every monomial ideal, and
 thus everywhere.  Every irreducible component must contain a monomial
 ideal, so Theorem \ref{monocase} says that every irreducible
 component must intersect the coherent component.  Since we just
 showed that $H_{\mathcal L}$ is smooth, we conclude that it is also
 irreducible.
\end{proof}

\noindent {\bf Acknowledgements.} We wish to thank Alastair Craw for 
details and references on $G$-Hilb.

\bibliographystyle{plain}
\bibliography{mt02} 

\begin{thebibliography}{10}

\bibitem{Arn}
V.I. Arnold.
\newblock A-graded algebras and continued fractions.
\newblock {\em Communications in Pure and Applied Mathematics}, 42:993--1000,
  1989.

\bibitem{BaMo}
D.~Bayer and I.~Morrison.
\newblock Gr{\"o}bner bases and geometric invariant theory {I}.
\newblock {\em Journal of Symbolic Computation}, 6:209--217, 1998.

\bibitem{BKR}
T.~Bridgeland, A.~King, and M.~Reid.
\newblock The {M}c{K}ay correspondence as an equivalence of derived categories.
\newblock {\em J. Amer. Math. Soc.}, 14(3):535--554 (electronic), 2001.

\bibitem{CrawReid}
A.~Craw and M.~Reid.
\newblock How to calculate {$A$}-{H}ilb {${\mathbb C}^3$}.
\newblock {math.AG/9909085}.

\bibitem{DHSS}
J.~Deloera, S.~Ho\c{s}ten, F.~Santos, and B.~Sturmfels.
\newblock The polytope of all triangulations of a point configuration.
\newblock {\em Documenta Mathematica}, 1:103--119, 1996.

\bibitem{GaPe}
V.~Gasharov and I.~Peeva.
\newblock Deformations of codimension 2 toric varieties.
\newblock {\em Compositio Mathematica}, 123:225--241, 2000.

\bibitem{HaSt}
M.~Haiman and B.~Sturmfels.
\newblock Multigraded {H}ilbert schemes.
\newblock math.AG/0201271.

\bibitem{HT1}
S.~Ho\c{s}ten and R.R. Thomas.
\newblock Standard pairs and group relaxations in integer programming.
\newblock {\em Journal of Pure and Applied Algebra}, 139:133--157, 1999.

\bibitem{Kid}
R.~Kidoh.
\newblock {H}ilbert schemes and cyclic quotient singularities.
\newblock {\em Hokkaido Mathematical Journal}, 30:91--103, 2001.

\bibitem{Lee}
C.~Lee.
\newblock Regular triangulations of convex polytopes.
\newblock In P.~Gritzmann and B.~Sturmfels, editors, {\em Applied Geometry and
  Discrete Mathematics - The Victor Klee Festschrift}, volume~4, pages
  443--456. AMS, Dimacs Series, Providence, R.I., 1991.

\bibitem{Mac}
D.~Maclagan.
\newblock Antichains of monomial ideals are finite.
\newblock {\em Proceedings of the AMS}, 129:1609--1615, 2001.

\bibitem{MT}
D.~Maclagan and R.~Thomas.
\newblock Combinatorics of the toric {H}ilbert scheme.
\newblock {\em Discrete and Computational Geometry}, 27:249--264, 2002.

\bibitem{MR}
T.~Mora and L.~Robbiano.
\newblock The {G}r{\"o}bner fan of an ideal.
\newblock {\em Journal of Symbolic Computation}, 6:183--208, 1988.

\bibitem{Oda}
T.~Oda.
\newblock {\em Convex Bodies and Algebraic Geometry}.
\newblock Springer-Verlag, Berlin, 1985.

\bibitem{PeSt1}
I.~Peeva and M.~Stillman.
\newblock Toric {H}ilbert schemes.
\newblock {\em Duke Math. J.}, 111(3):419--449, 2002.

\bibitem{SST}
M.~Saito, B.~Sturmfels, and N.~Takayama.
\newblock {\em Gr{\"o}bner Deformations of Hypergeometric Differential
  Equations}, volume~6.
\newblock Algorithms and Computation in Mathematics, Springer-Verlag, New-York,
  1999.

\bibitem{Stu94}
B.~Sturmfels.
\newblock The geometry of ${A}$-graded algebras.
\newblock 1994.
\newblock math.AG/9410032.

\bibitem{GBCP}
B.~Sturmfels.
\newblock {\em Gr{\"o}bner Bases and Convex Polytopes}.
\newblock American Mathematical Society, Providence, RI, 1995.

\end{thebibliography}

\end{document}